\documentclass[reqno]{amsart}
\usepackage{amsmath}
\usepackage{amscd}
\usepackage{graphics}
\usepackage{latexsym}

\theoremstyle{plain}
\newtheorem{theorem}{Theorem}[section]
\newtheorem{thm}[theorem]{Theorem}
\newtheorem{cor}[theorem]{Corollary}
\newtheorem{lem}[theorem]{Lemma}
\newtheorem{prop}[theorem]{Proposition}

\theoremstyle{definition}

\theoremstyle{remark}

\newcommand{\CC}{\mathbb{C}}

\newcommand{\PP}{\mathbb{P}}

\newcommand{\Kbm}[2]{\overline{\mathcal{M}}_{#1}(#2)}

\newcommand{\lt}{\left}
\newcommand{\rt}{\right}

\newcommand{\mc}{\mathcal}

\newcommand{\Metc}{M^{d,w}}
\newcommand{\Hetc}{H^{d,w}}
\newcommand{\Betc}{(B,\Sigma,b_0)}

\begin{document}

\title{A note on Hurwitz schemes of covers of a positive genus curve}
\author[Harris]{Joe Harris} 
\address{Department of Mathematics \\ 
  Harvard University \\ Cambridge MA 02138}
\email{harris@math.harvard.edu} 
\author[Graber]{Tom Graber }
\address{Department of Mathematics \\ 
  Harvard University \\ Cambridge MA 02138}
\email{graber@math.harvard.edu}
\author[Starr]{Jason Starr}
\address{Department of Mathematics \\ 
  Massachusetts Institute of Technology \\ Cambridge MA 02139}
\email{jstarr@math.mit.edu} \date{\today}

\begin{abstract}
  Let $B$ be a smooth, connected, projective complex curve of genus
  $h$.  For $w\geq 2d$ we prove the irreducibility of the Hurwitz
  stack ${\mc H}_{S_d}^{d,w}(B)$ parametrizing degree $d$ covers
  of $B$ simply-branched over $w$ points, and with monodromy group $S_d$.
\end{abstract}

\maketitle

\section{Introduction}~\label{sec-intro}
Suppose that $B$ is a smooth, connected, projective complex curve of
genus $h$.  Let $d >0$ and $w\geq 0$ be integers such that $g :=
d(h-1) + \frac{w}{2} + 1$ is a nonnegative integer 
(in particular $w$ is even).
We define ${\mc H}^{d,w}(B)$ to be the open substack of the Kontsevich
moduli stack $\Kbm{g,0}{B,d}$ parametrizing stable maps
$f:X\rightarrow B$ such that $X$ is smooth and $f$ is finite with only
simple branching.  Let $\text{br}(f) \subset B$ denote the branch
divisor of $f$.
If we choose a basepoint $b_0\in B-\text{br}(f)$ and an
identification $\phi:f^{-1}(b_0) \rightarrow \{1,\dots, d\}$, there is
an induced monodromy homomorphism $\tilde{\phi}:\pi_1(B-\text{br}(f),
b_0) \rightarrow S_d$ which associates to any loop
$\gamma:[0,1]\rightarrow B$ with $\gamma(0)=\gamma(1)=b_0$, the
permutation of $f^{-1}(b_0)$ determined by analytic continuation along
$\gamma$.  In particular, the subgroup
$\text{image}(\tilde{\phi})\subset S_d$ is well-defined up to
conjugation independently of $\phi$.  The corresponding conjugacy
class of subgroups determines a locally constant function on ${\mc
  H}^{d,w}(B)$.  Given a subgroup $G\subset S_d$, we define 
${\mc H}_G^{d,w}(B)$ to be the open and closed substack of 
${\mc H}^{d,w}(B)$ parametrizing stable maps $f:X\rightarrow B$ whose
corresponding monodromy group is conjugate to $G$.  We are
particularly interested in ${\mc H}^{d,w}_{S_d}(B)$, the stack parametrizing
Hurwitz covers of $B$ with full monodromy group.

\

\begin{thm}~\label{thm-thm1}
If $w\geq 2d$, then ${\mc H}^{d,w}_{S_d}(B)$ is a connected, smooth,
finite-type Deligne-Mumford stack over $\CC$.
\end{thm}

The fact that ${\mc H}^{d,w}_{S_d}(B)$ is a finite-type Deligne-Mumford
stack follows from the fact that $\Kbm{g,0}{B,d}$ is a finite-type
Deligne-Mumford stack.  The fact that ${\mc H}^{d,w}_{S_d}(B)$ is smooth
follows from a trivial deformation theory computation.  So the content
of theorem~\ref{thm-thm1} is that ${\mc H}^{d,w}_{S_d}(B)$ is connected.

\

This is a classical fact when $h=0$, i.e. for branched covers of
$\PP^1$, (c.f. \cite{Clebsch}, \cite{Hurw}, and for a modern account
\cite[prop. 1.5]{F69}).
This fact is well-known to experts, but there seems to be no
reference.  We used theorem~\ref{thm-thm1} in our paper~\cite{GHS}, and so we
present a proof below.  We wish to thank Ravi Vakil for useful
discussions.

\section{Setup}~\label{sec-setup}

Our eventual goal is to prove theorem~\ref{thm-thm1},
but for most of this paper, we shall work
with schemes which admit \'etale maps to ${\mc H}^{d,w}_{S_d}(B)$.
Suppose
$\Sigma\subset B$ is a finite subset, and suppose $b_0\in \Sigma$ is a
point.  
We define $\Metc\Betc$ to be the fine moduli scheme
parametrizing pairs $(f:X\rightarrow B,\phi)$ where $f:X\rightarrow B$
is a stable map in $\Kbm{g,0}{B,d}$ and where
$\phi:f^{-1}(b_0)\rightarrow \{1,\dots,d\}$ are such that
\begin{enumerate}
\item $f$ is finite,
\item $f$ is unramified over $\Sigma$, and
\item $\phi$ is a bijection.
\end{enumerate}

\

Using known results on the Kontsevich moduli space $\Kbm{g,0}{B,d}$,
it is easy to show that $\Metc\Betc$ is a nonempty, smooth,
quasi-projective scheme of dimension $w$.
By ~\cite{FaP}, there is a branch morphism
$\text{br}:\Metc\Betc\rightarrow (B-\Sigma)_w$ where 
$(B-\Sigma)_w$ is the $w$th symmetric power parametrizing
effective degree $w$ divisors on $B-\Sigma$.  It is clear that
$\text{br}$ is quasi-finite, and thus $\text{br}:\Metc\Betc\rightarrow
(B-\Sigma)_w$ is dominant.  We denote by $(B-\Sigma)_w^o\subset
(B-\Sigma)_w$ the Zariski open subset parametrizing reduced
effective divisors of degree $w$ in $B-\Sigma$.  We define
$\Hetc\Betc\subset \Metc\Betc$ to be the preimage under $\text{br}$ of
$(B-\Sigma)_w^o$.  

\

For each
pair $(f:X\rightarrow B,\phi)$ in $\Hetc\Betc$ with branch divisor
$\text{br}(f)$, there is an induced monodromy homomorphism
$\tilde\phi:\pi_1(B-\text{br}(f),b_0)\rightarrow S_d$ where $S_d$ is
the symmetric group of permutations of $\{1,\dots,d\}$.  The image of
$\tilde\phi$ determines a locally constant function on $\Hetc\Betc$.
Because $\Metc\Betc$ is smooth and $\Hetc\Betc$ is dense in
$\Metc\Betc$, this locally constant function extends to all of
$\Metc\Betc$.  Given a subgroup $G\subset S_d$ we define
$\Metc_G\Betc$ (resp. $\Hetc_G\Betc$) 
to be the open and closed subscheme of $\Metc\Betc$ (resp. $\Hetc\Betc$)
on which the image of $\tilde\phi$ equals $G$.

\

Let $F:{\mc X}(w)\rightarrow \Hetc_{S_d}\Betc \times B$ be the
pullback of the universal stable map, i.e. ${\mc X}(w)$ parametrizes
data $(f:X\rightarrow B,\phi,x)$ 
where $(f:X\rightarrow B,\phi)\in \Hetc_{S_d}\Betc$ and $x\in X$, and
$F(f:X\rightarrow B,\phi,x) = (f:X\rightarrow B,\phi,f(x))$.  
We denote by $U\subset {\mc X}(w)\times_{F,F} {\mc X}(w)$ the open
subscheme of the fiber product of
${\mc X}(w)$ with itself over $\Hetc_{S_d}\Betc \times B$ 
parametrizing data $(f:X\rightarrow B,\phi,x_1,x_2)$ such that $x_1\neq
x_2$, and such that $f(x_1)=f(x_2)$ is neither in $S$ nor equal to any
branch point of $f$.  We define ${\mc X}_2(w)$ to be the quotient of $U$ by
the obvious involution $(f:X\rightarrow B,\phi,x_1,x_2)\sim
(f:X\rightarrow B,\phi,x_2,x_1)$.  We denote by ${\mc X}_2^e(w)$ the open
subscheme of the $e$-fold fiber product of ${\mc X}_2(w)$ with itself over
$\Hetc_{S_d}\Betc$ parametrizing data $(f:X\rightarrow
B,\phi,\{x^1_1,x^2_2\}, \dots,\{x^e_1,x^e_2\})$ such that
$f(x^1_1),\dots,f(x^e_1)$ are all distinct points in $B-\Sigma$.  
Notice that the projection ${\mc X}_2^e(w)\rightarrow
\Hetc_{S_d}\Betc$ is flat.  The condition that the image of
$\tilde\phi$ be all of $S_d$, and therefore doubly-transitive, implies
that ${\mc X}_2(w)\rightarrow \Hetc_{S_d}\Betc$ has irreducible
fibers.  Therefore also ${\mc X}_2^e(w)\rightarrow \Hetc_{S_d}\Betc$
has irreducible fibers.

\

For each $(f:X\rightarrow
B,\phi,\{x^1_1,x^1_2\},\dots,\{x^e_1,x^e_2\})$ in ${\mc X}_2^e(w)$ we
can associate a pair 
$(f_a:X_a\rightarrow B,\phi_a)$ in $H^{d,w+2e}_{S_d}\Betc$ as follows: 
\begin{enumerate}
\item  We define $X_a$ to be 
the $e$-nodal curve whose normalization is of the form $u:X\rightarrow
X_a$ such
that $u(x^i_1)=u(x^i_2)$ for each $i=1,\dots,e$,
\item we define $f_a:X_a\rightarrow B$ to be the unique
morphism such that $f=f_a\circ u$, and
\item we define $\phi_a$ to be the unique map
such that $\phi=\phi_a\circ u$.  
\end{enumerate}
This association defines a regular
morphism $G_{w,e}:{\mc X}_2^e(w)\rightarrow M^{d,w+2e}_{S_d}\Betc$.  

\

\begin{figure}
\begin{picture}(330,170)
\put(10,10){\makebox(330,140){\includegraphics{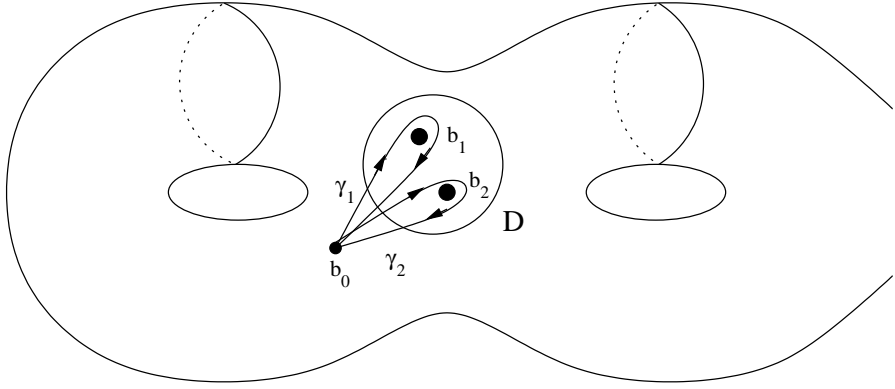}}}
\end{picture}
\caption{Adjacent branch points with equal monodromy}
\label{fig-F1}
\end{figure}

\

We give a topological application of the morphism $G_{w,e}$.  Suppose
we have a pair $(f:X\rightarrow B,\phi)$ in $H^{d,w+2}_{S_d}\Betc$ 
and $D\subset B$ is a closed disk which is disjoint from $\Sigma$,
such that $D\cap \text{br}(f)$ consists of two branch points $b_1,b_2$
which are contained in the interior of $D$.  Define $U=B-D$ and
suppose that 
$\tilde\phi:\pi_1(U-\text{br}(f),b_0)\rightarrow S_d$ is surjective.
Suppose moreover that $f$ is trivial over the boundary
$\partial D$ of $D$, i.e. $f^{-1}(\partial D)$ consists of $d$ disjoint
circles each of which maps homeomorphically to $\partial D$.  Choose
simple closed loops $\gamma_1,\gamma_1$ around $b_1$ and $b_2$ 
as displayed in Figure~\ref{fig-F1}.
Then $\tilde\phi(\gamma_1)$ and $\tilde\phi(\gamma_2)$ both equal the
same transposition $\tau = (j,k)$.  

\

\begin{lem}~\label{lem-split}  With the notations and assumptions in
  the last paragraph, there exists a pair $(f_a:X_a\rightarrow
  B,\phi_a)$ in $H^{d,w}_{S_d}\Betc$, a datum $(f_a:X_a \rightarrow B,
  \phi_a, \{x_j, x_k\})$ in ${\mc X}_2(w)$, and an
  analytic isomorphism $h:f^{-1}(U) \rightarrow f_a^{-1}(U)$ such that
\begin{enumerate}
\item $f|_{f^{-1}(U)} = (f_a)|_{f_a^{-1}(U)}\circ h$ and $\phi =
  \phi_a\circ h$, and
\item the image by $G_{w,1}$ of $(f_a:X_a\rightarrow B, \phi_a,
  \{x_j,x_k\})$ lies in the same connected component of
  $M^{d,w+2}_{S_d}\Betc$ as $(f:X\rightarrow B,\phi)$.
\end{enumerate}
\end{lem}

\begin{proof}
We may choose an analytic isomorphism of the disk $D\subset B$ with
the unit disk $\Delta\subset \CC$ such that $b_1$ and $b_2$ map to the
two roots of $x^2=t_0$ for some $t_0\in \Delta-\{0\}$.  Let $x$ be the
coordinate on $\Delta$.
Consider the map
$f^{-1}(D)\rightarrow D$.  For each $i\neq j,k$ the connected
component of $f^{-1}(D)$ corresponding to $i$ maps isomorphically to
$D$.  The connected component of $f^{-1}(D)$ corresponding to $j$ and
$k$ is identified with the covering $C_{t_0}$ 
of $\Delta$ given by $C_{t_0} = \{(x,y)\in
\CC^2: x\in \Delta, y^2-(x^2-t_0)=0\}$.  For $t\in \Delta$, consider
the family of covers $C_t=\{(x,y)\in \CC^2: x\in \Delta, y^2
-(x^2-t)=0\}$.  By the Riemann existence theorem, for each $t\in
\Delta$ there is a pair $(f_t:X_t\rightarrow B,\phi_t)$ in
$M^{d,w+2}\Betc$ such that the
restriction of $f_t$ to $f_t^{-1}(U)$ is identified with the
restriction of $f$ to $f^{-1}(U)$ and such that the restriction
of $f_t$ to $D$ consists of $d-2$ copies of $D$ mapping isomorphically
to $D$ (one copy for each $i\neq j,k$), and the connected component
corresponding to $j$ and $k$ is identified with $C_t\rightarrow
\Delta$.  We will see that $(f_0:X_0\rightarrow B,\phi_0)$ is in the
image of $G_{w,1}: {\mc X}_2(w)\rightarrow M^{d,w+2}\Betc$.  

\

Define $u:X_a\rightarrow X_0$ to be the normalization and define
$\{x_j,x_k\}$ to be the preimage of the node $x_0\in X_0$.  We define
$f_a:X_a\rightarrow B$ to be $f_a=f_0\circ u$ and $\phi_a=\phi_0\circ
u$.  Notice that $f_a:X_a\rightarrow B$ is unbranched over $D$.
Define $x_j$ (resp. $x_k$) to be the preimage of $f_0(x_0)$
on the sheet of $f_a^{-1}(D)$ corresponding to $j\in\{1,\dots,d\}$
(resp. to $k\in \{1,\dots, d\}$).
We have an identification of $u^{-1}(f_0^{-1}(U))$ with $f_0^{-1}(U)$.
Therefore we have an identification $h:f^{-1}(U)\rightarrow
f_a^{-1}(U)$ commuting with $f, f_a$ and with $\phi, \phi_a$.  In
particular, we conclude that
$\tilde\phi_a:\pi_1(U-\text{br}(f_a),b_0)\rightarrow S_d$ is identified
with $\tilde\phi:\pi_1(U-\text{br}(f),b_0)\rightarrow S_d$ and so is
surjective.  So $(f_a:X_a\rightarrow B,\phi_a)$ is in
$H^{d,w}_{S_d}\Betc$.  Clearly $(f_a:X_a\rightarrow B,\phi_a,
\{x_j,x_k\})$ is in ${\mc X}_2(w)$ and, by construction, its image
under $G_{w,1}$ is 
$(f_0:X_0\rightarrow B,\phi_0)$.  Since $(f_0:X_0\rightarrow
B,\phi_0)$ is in the same connected component of
$H^{d,w+2}_{S_d}\Betc$ as $(f:X\rightarrow B,\phi)$, this proves the
lemma. 
\end{proof}

\

\begin{lem}~\label{lem-dbltran}  With the notations and assumptions in
  lemma~\ref{lem-split}, suppose given a transposition $(j_b,k_b)\in
  S_d$.
  Then there exists a pair $(f_b:X_b\rightarrow B,\phi_b)$ in
  $H^{d,w+2}_{S_d}\Betc$, and an analytic isomorphism $h:f^{-1}(U)
  \rightarrow f_b^{-1}(U)$ such that:
\begin{enumerate}
\item $\text{br}(f_b)=\text{br}(f)$,
\item $f|_{f^{-1}(U)} = f_b|_{f_b^{-1}(U)}\circ h$ and
  $\phi = \phi_b \circ h$, 
\item $\tilde{\phi_b}(\gamma_1)=\tilde{\phi_b}(\gamma_2) = (j_b,k_b)$, and
\item $(f_b:X_b\rightarrow B,\phi_b)$ is in the same connected
  component of $H^{d,w+1}_{S_d}\Betc$ as $(f:X\rightarrow B,\phi)$.
\end{enumerate}
\end{lem}

\begin{proof}
Let $(f_a:X_a\rightarrow B,\phi_a, \{x_j,x_k\})$,
$h_a:f^{-1}(U)\rightarrow f_a^{-1}(U)$ be as constructed in
the proof of lemma~\ref{lem-split}.  Define $b_1$ to be
$f_a(x_j)=f_a(x_k)$.  Define $(X_a)_2\rightarrow B- \lt(\Sigma \cup
\text{br}(f)\rt)$   
to be the fiber of
${\mc X}_2(w)\rightarrow \Hetc_{S_d}\Betc$ over $(f_a:X_a\rightarrow
B,\phi_a)$.  Notice that $(X_a)_2 \rightarrow B-\lt(\Sigma \cup
\text{br}(f)\rt)$ is an
unbranched covering space.  Define $x_{j_b}$ (resp. $x_{k_b}$) to be
the elements of $f_a^{-1}(b_1)$ which lie on the sheets of
$f_a^{-1}(D)$ corresponding to $j_b\in \{1,\dots, d\}$ (resp. $k_b\in
\{1,\dots, d\}$). 
Because $\tilde\phi_a:\pi_1(B-\text{br}(f_a),b_1)\rightarrow S_d$ is
surjective, in particular it is doubly transitive.  
Therefore $(X_a)_2$ is irreducible and
$(f_a:X_a \rightarrow B,
\phi_a, \{x_{j_b}, x_{k_b}\})$ is in the same connected component of
${\mc X}_2(w)$ as $(f_a:X_a \rightarrow B, \phi_a, \{x_j, x_k\})$.  So
the image $G_{w,1}(f_a:X_a \rightarrow B,
\phi_a, \{x_{j_b}, x_{k_b}\})$ is in the same connected component of
$M^{d,w+2}_{S_d}\Betc$ as $(f:X \rightarrow B, \phi)$.  

\

Consider the same family of covers $C_t\rightarrow \Delta$
as in the proof of lemma~\ref{lem-split}, with the roles of $j,k$
replaced by $j_b,k_b$. By the Riemann existence theorem 
there exists a pair $(f_b:X_b \rightarrow B, \phi_b)$ and an isomorphism
$h_b:f_a^{-1}(U)\rightarrow f_b^{-1}(U)$ commuting with $f_a, f_b$ and
$\phi_a, \phi_b$
such that $f_b^{-1}(D)\rightarrow D$ is identified with the covering
$C_1\rightarrow \Delta$.  We define $h:f^{-1}(U)\rightarrow
f_b^{-1}(U)$ to be $h_b\circ h_a$.  Then $(f_b:X_b \rightarrow B,
\phi_b)$ and $h$ satisfy items $(1)$, $(2)$, and $(3)$ of the lemma.
Moreover $(f_b:X_b\rightarrow B, \phi_b)$ is in the same connected
component of $M^{d,w+2}_{S_d}\Betc$ as $G_{w,2}(f_a:X_a \rightarrow B,
\phi_a, \{x_{j_b},x_{k_b}\} )$.  Thus $(f_b:X_b \rightarrow B,
\phi_b)$ is in the same connected component of $M^{d,w+2}_{S_d}\Betc$
as $(f:X \rightarrow B, \phi)$.  Since both pairs are in
$H^{d,w+2}_{S_d}\Betc$ and since $M^{d,w+2}_{S_d}\Betc$ is smooth, we
conclude that both pairs are in the same connected component of
$H^{d,w+2}_{S_d}\Betc$.  
\end{proof}

\section{Branching monodromy}~\label{sec-brmon}

Fix a closed disk $D\subset B$ disjoint from $\Sigma$.  Fix a path from
$b_0$ to the boundary $\partial D$.  Denote by $U$ the open subset
$B-D\subset B$.
In most of this section we will restrict our attention to the
analytic open subset $V\subset \Metc_{S_d}\Betc$ parametrizing
$(f:X\rightarrow B,\phi)$ such that $\text{br}(f)$ is contained in the
interior of $D$.
By assumption, the monodromy group of $f$, i.e. the image of
$\tilde\phi$, is all of $S_d$.  For each connected component of 
$\Hetc_{S_d}\Betc \cap V$, the function which associates to each
$(f:X\rightarrow B,\phi)$ the image of
$\tilde\phi:\pi_1(D-\text{br}(f),b_0)\rightarrow S_d$ is constant.
We call this subgroup the \emph{branching monodromy group} of $f$ (of
course it depends on the choice of $D$ and the path from $b_0$ to
$\partial D$).

\

Since $\Hetc_{S_d}\Betc \cap V$ is the complement of proper analytic
subvarieties of the complex manifold $V$, each connected component of
$V$ is the closure of a unique connected component of 
$\Hetc_{S_d}\Betc \cap V$.  For each subgroup $G\subset D$
let us denote by $V_G\subset V$ the open and closed submanifold on
which the image of $\tilde\phi:\pi_1(D-\text{br}(f),b_0)\rightarrow
S_d)$ equals $G$.  The goal of this section is to prove that when
$w\geq 2d$, every connected component of $\Hetc_{S_d}\Betc$ has
nonempty intersection with $V_{S_d}$, i.e. there is a pair
$(f:X\rightarrow B,\phi)$ in this connected component and in $V$ which
has branching monodromy group equal to $S_d$.  

\

Suppose that $(f:X\rightarrow B,\phi)\in \Hetc_{S_d}\Betc \cap V_G$. 
Because $G$ is generated by
transpositions, there is a partition $(A_1,\dots,A_r)$ of
$\{1,\dots,d\}$ such that $G=S_{A_1}\times \dots \times S_{A_r}$ where
$S_{A_m}\subset S_d$ consists of those permutations which act as the
identity on each subset $A_n\subset \{1,\dots,d\}$ for which $n\neq
m$.  In other words, $G$ is the subgroup of permutations which
stabilize each subset $A_m\subset \{1,\dots,d\}$.

\

\begin{figure}
\begin{picture}(330,240)
\put(0,10){\makebox(330,210){\includegraphics{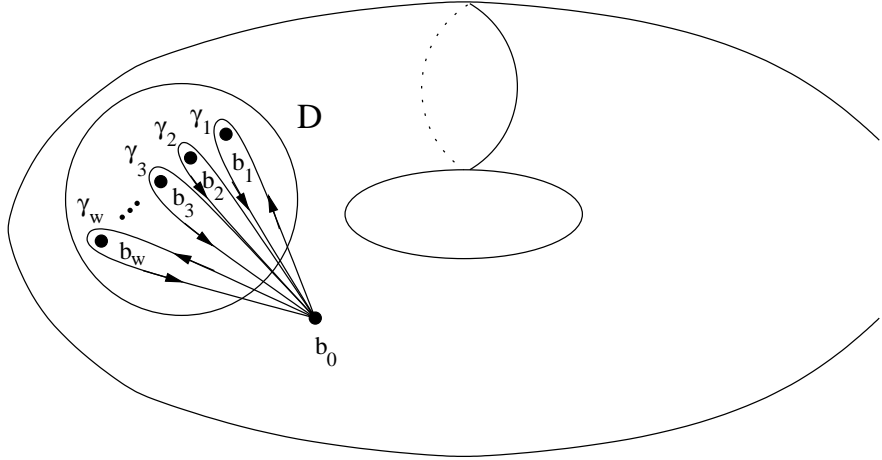}}}
\end{picture}
\caption{Branch points contained in the disk $D$}
\label{fig-P1}
\end{figure}

\

Choose a system of loops $\gamma_1,\dots, \gamma_w$
as in Figure~\ref{fig-P1}.  Denote by $\tau_i$ the transposition
$\tilde{\phi}(\gamma_i)$.  Then each $\tau_i$ lies in one of the
subgroups $S_{A_{m(j)}}$.  

\

\begin{figure}
\begin{picture}(330,280)
\put(10,10){\makebox(330,250){\includegraphics{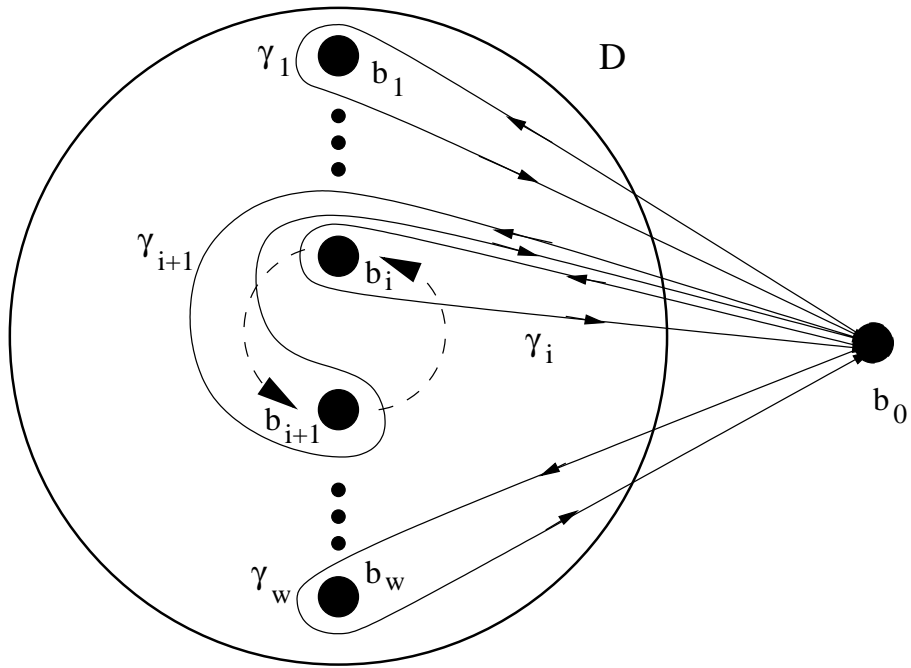}}}
\end{picture}
\caption{Braid move exchanging two branch points}
\label{fig-P1a}
\end{figure}
\

\

\begin{figure}
\begin{picture}(330,280)
\put(10,10){\makebox(330,250){\includegraphics{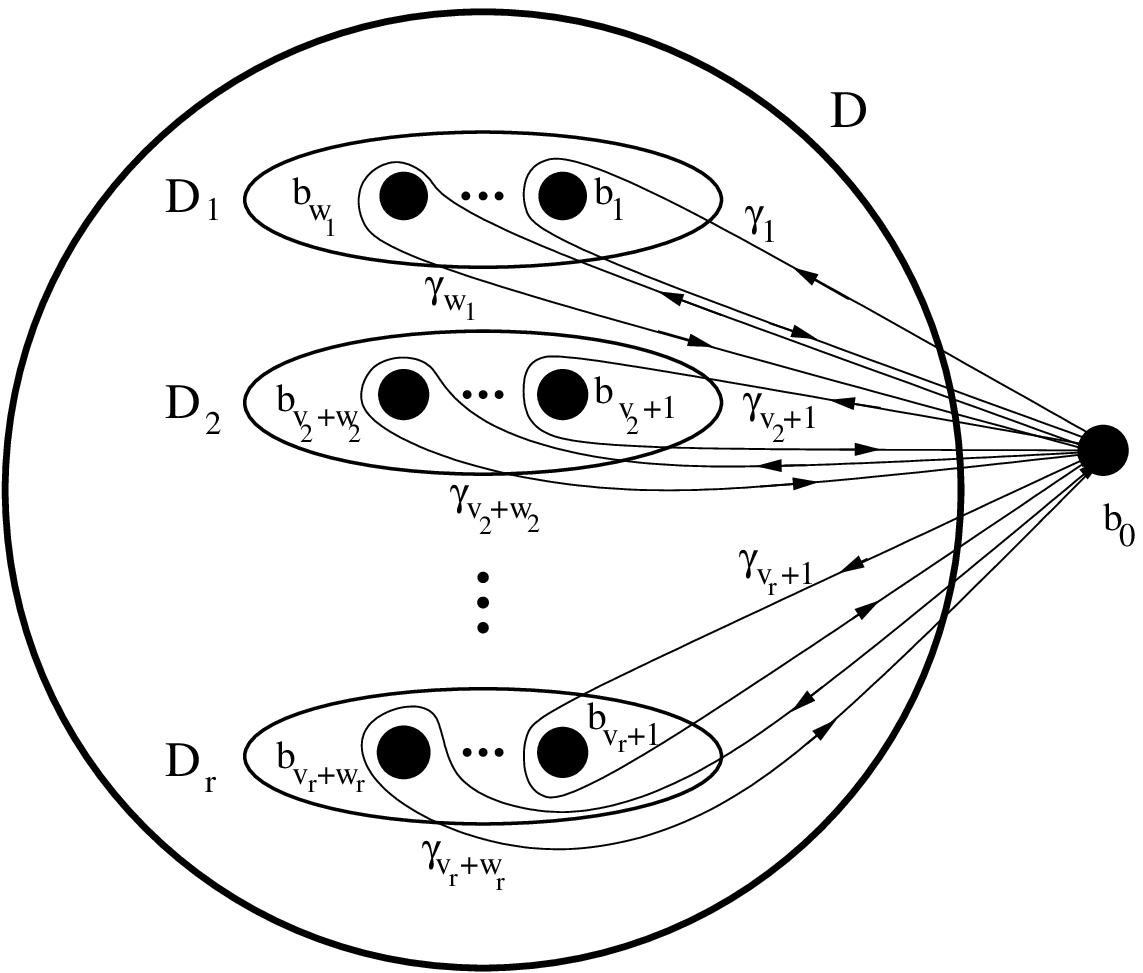}}}
\end{picture}
\caption{Branch points in standard position}
\label{fig-P3}
\end{figure}
\ 

Suppose that $\gamma_i$ and $\gamma_{i+1}$
are adjacent loops such that $\tau_i$ lies in $S_{A_m}$ and
$\tau_{i+1}$ lies in $S_{A_n}$ with $m\neq n$.  Consider the element
$\sigma_i$ 
of the braid group which interchanges the branch points
$b_i$ and $b_{i+1}$ as shown in Figure~\ref{fig-P1a}.  The result is to replace
$\tau_i$ by $\tau_{i+1}$ and to replace $\tau_{i+1}$ by
$\tau_{i+1}\tau_i\tau_{i+1}$.  Since $A_m$ and $A_n$ are disjoint, we
have $\tau_{i+1}\tau_i\tau_{i+1}=\tau_i$.  In other words, the result
is to interchange $\tau_i$ and $\tau_{i+1}$.  Note that this
operation does not change $G\subset S_d$. By repeating this
process, we may arrange that there are integers $w_0=0,w_1,\dots,w_r$
with the following property: for $m=1,\dots, r$, denote $v_m = w_0 +
\dots + w_{m-1}$;
then for each $m=1,\dots,r$, each transpositions $\tau_i$ with $v_m+1
\leq i \leq v_{m+1}$ is in
$S_{A_m}$. 
Notice that since these transpositions generate $S_{A_m}$, we have
$w_m \geq \#A_m - 1$.  
Stated more precisely, we have proved that given
a pair $(f:X\rightarrow B,\phi)$ in $\Hetc_{S_d} \Betc \cap V_G$, 
in the same connected
component of $\Hetc_{S_d} \Betc \cap V_G$ 
there is a pair $(f_a:X_a \rightarrow B, \phi_a)$ and an isomorphism
$h:f^{-1}(U) \rightarrow (f_a)^{-1}(U)$
such that:
\begin{enumerate}
\item $\text{br}(f_a) = \text{br}(f)$,
\item $f|_{f^{-1}(U)} =
(f_a)|_{f_a^{-1}(U)} \circ h$ and
$\phi=\phi_a\circ h$, and 
\item the transpositions $\tau_i=\tilde{\psi}(\gamma_i)$ satisfy
  $\gamma_i\in S_{A_m}$ for $v_m+1\leq i \leq
  v_{m+1}$. 
\end{enumerate}
We say that a pair $(f_a:X_a \rightarrow B, \phi_a)$ 
satisfying item $(3)$ is in \textit{standard position}.  
For each $m=1,\dots,r$, choose a subdisk $D_m\subset D$
as in Figure~\ref{fig-P3} which contains the loops $\gamma_i$ for
$(w_0+\dots+w_{m-1})+1\leq i \leq w_0+\dots+w_{m-1}+w_m$.
Note that any braid move in $D_m$ has no effect on the branch points belonging
to $D_n$ with $n\neq m$.  

\

\begin{prop}~\label{prop-split}  Suppose that $(f:X\rightarrow
  B,\phi)$ in $\Hetc_{S_d}\Betc \cap V_G$ is in standard position.  
  Suppose $w_m\geq
  2\#A_m$.  Then there are braid moves in $D_m$ transforming
  $(\tau_{v_m+1},\dots, \tau_{v_m+w_m}$ into
  $(\tau'_1,\dots,\tau'_{w_m-2},\tau,\tau)$ such that
  $\tau'_1,\dots,\tau'_{w_m-2}$ generate $S_{A_m}$. 
\end{prop}

\begin{proof}
Define $g=\tau_{v_m+1}\cdot \dots \cdot \tau_{v_{m+1}}$.  
By ~\cite[theorem 1]{Klu}, the braid group of $D_m$ acts transitively
on the set 
\begin{eqnarray*}
O_g:=\{(\tau_1,\dots,\tau_{w_m})\in S_{A_m} | \text{ each }
\tau_i \text{ a transposition }, \\
\langle \tau_1,\dots,\tau_{w_m}
\rangle = S_{A_m}, \tau_1\cdot\dots \cdot \tau_{w_m}=g \}.
\end{eqnarray*}
Thus it suffices to find $(\tau'_1,\dots,\tau'_{w_m-2},\tau,\tau)$ as
above which lies in $O_g$.  

\

Suppose that $g$ has cycle type $(\lambda_1,\dots,\lambda_s)$ 
for some partition
$\lambda= (\lambda_1\geq \lambda_2 \geq \dots \geq \lambda_s)$ of
$\#A_m$. 
Define $\lambda_0=0$ and for each $k=1,\dots,s$, define $\mu_k =
\lambda_0+\dots + \lambda_{k-1}$.  Then we may order the elements of
$A_m$ so that $g$ is the permutation 
\begin{equation}
g = \lt(\mu_1+1,\dots,\mu_1+\lambda_1\rt)
\lt(\mu_2+1,\dots,\mu_2+\lambda_2\rt) \dots 
\lt(\mu_s+1,\dots,\mu_s+\lambda_2\rt).
\end{equation}
Of course this ordering has nothing to do with the ordering induced by
$\phi$.  

\

For each $k=1,\dots,s$, consider the ordered sequence of
transpositions, which is defined to be empty if $\lambda_k=1$, and for
$k>1$ is
defined to be 
\begin{equation}
I_k = \lt( (\mu_k+1,\mu_k+2), (\mu_k+2,\mu_k+3), \dots ,
(\mu_k+\lambda_k-1,\mu_k+\lambda_k) \rt).
\end{equation}
Thus $I_k$ contains $\lambda_k-1$ transpositions.  Next consider the
sequence of transpositions
\begin{equation}
I_{s+1} = \lt( (\mu_1,\mu_2), (\mu_1,\mu_2), (\mu_2,\mu_3), (\mu_2,\mu_3),
\dots, (\mu_{s-1},\mu_s),(\mu_{s-1},\mu_s) \rt).
\end{equation}
The concatenated sequence $I = I_1\cup \dots \cup I_s \cup I_{s+1}$ has
length $L:=\sum_k (\lambda_k -1 ) + 2(s-1) = \#A_m + s - 2 \leq 2\#A_m -2$.
The product of these transpositions is $g$, and these transpositions
generate $S_{A_m}$.  Since the sign of $g$ is both $(-1)^{w_m}$ and
$(-1)^L$, we have that $w_m-L$ is divisible by $2$.  And the
assumption that $w_m \geq 2\#A_m$, implies that $w_m - L \geq 2$.
If we choose any transposition $\tau\in S_{A_m}$,
and let $J$ be the constant sequence of length $w_m-L$,
$J=(\tau,\tau,\dots, \tau)$, then we have that the concatenated
sequence $I\cup J$ is an element of $O_g$ satisfying the hypotheses of
the proposition.
\end{proof}

\

\begin{cor}~\label{cor-cor1}
Given $w'\geq w$ with $w'\geq 2d, w\geq 2d-2$, set $e=\frac{w'-w}{2}$.  
Suppose given a pair
$(f:X\rightarrow B,\phi)\in H^{d,w'}_{S_d}\Betc \cap V_G$.  Then  
there is a pair $(f_a:X_a \rightarrow B, \phi_a)\in H^{d,w}_{S_d}\Betc
\cap V_G$, an isomorphism
$h:f^{-1}(U)\rightarrow (f_a)^{-1}(U)$ and a datum 
$(f_a:X_a \rightarrow B, \phi_a, \{x^1_1,x^1_2\},
  \dots, \{x^e_1,x^e_2\})$ in ${\mc X}_2^e(w)$ 
such that
\begin{enumerate}
\item $f_a(x^i_j)\in D$ for $i=1,\dots, e$ and for $j=1,2$,
\item $f|_{f^{-1}(U)}=(f_a)|_{f_a^{-1}(U)}\circ h$ and $\phi=\psi\circ h$, and
\item the image of 
  $(f_a:X_a \rightarrow B, \phi_a, \{x^1_1,x^1_2\},
  \dots, \{x^e_1,x^e_2\})$ 
  under $G_{w,e}$ is contained in the same connected component
  of $H^{d,w'}_{S_d}\Betc$ as $(f:X\rightarrow B,\phi)$.  
\end{enumerate}
\end{cor}

\begin{proof}  We prove this by induction on $w'-w$.  For $w'=w$, there
  is nothing to prove.  Suppose $w'-w > 0$ and suppose the proposition
  has been proved for all smaller values of $w'-w$.  We note by
  proposition~\ref{prop-split} that there is a map
  $(f_c:X_c\rightarrow B, \phi_c)$ and $h_c:f^{-1}(U)\rightarrow
  f_c^{-1}(U)$ satisfying the conditions of that proposition.  If we
  define $D'$ to be a small disk containing the branch points of $f_c$
  corresponding to the transposition $\tau$, 
  then $(f_c:X_c\rightarrow B, \phi_c)$ and $D'$ satisfy the
  hypothesis of lemma~\ref{lem-split}.  By that lemma,
  there is a datum $(f_b:X_b \rightarrow B,\phi_b,\{x_j,x_k\}) \in
  {\mc X}_2(w'-2)$ and an isomorphism $h_b:f_a^{-1}(B-D')\rightarrow
  f_b^{-1}(B-D')$ satisfying the conditions of that lemma.  

\

If $w' = w+2$, we are done by taking $(f_a:X_a \rightarrow B, \phi_a,
  \{x^1_1, x^1_2\} ) = (f_b: X_b \rightarrow B, \phi_b, \{x_j, x_k\})$
  and taking $h = h_b \circ h_c$.  Therefore suppose that $w' > w+2$.  
  Now $(f_b:X_b \rightarrow B, \phi_b)$ is in $H^{d,w'-2}_{S_d}\Betc
  \cap V_G$.
  Since $(w'-2)-w < w'-w$, by the induction assumption
  there exists a datum $(f_a:X_a \rightarrow B, \phi_a,
  \{x^1_1,x^1_2\}, \dots, \{x^{e-1}_1,x^{e-1}_2\})$ in ${\mc
  X}_2^{e-1}(w)$ and
  $h_a:f_b^{-1}(U)\rightarrow f_a^{-1}(U)$ satisfying the conditions
  of our corollary.  Up to deforming this datum slightly, we may
  suppose that the isomorphism $h_a$ extends to a larger open set
  which contains $x_j, x_k\in X_b$, and, defining $x^e_1 =
  h_a^{-1}(x_j)$ and $x^e_2 = h_a^{-1}(x_k)$, the datum 
  $(f_a: X_a \rightarrow B,
  \phi_a, \{ x^1_1, x^1_2 \} ,\dots, \{x^{e-1}_1, x^{e-1}_2 \},
  \{ x^e_1, x^e_2 \} )$ is in ${\mc X}_2^e(w)$.  We
  define $h=h_a\circ h_b \circ h_c$.  The image of this datum
  under $G_{w,e}$ is contained in the same connected component as the
  image of $(f_b:X_b \rightarrow B, \phi_b, \{x_j, x_k\})$ under
  $G_{w'-2,1}$.  So the corollary is proved by induction.
\end{proof}

\

\begin{cor}~\label{cor-maincor}
If $w\geq 2d$, then for any pair $(f:X\rightarrow B,\phi)\in
\Hetc_{S_d}\Betc$ in $V$, there is a pair $(f_a:X_a \rightarrow B, \phi_a)
\in
\Hetc_{S_d}\Betc \cap V_{S_d}$ and an isomorphism $h:f^{-1}(U)\rightarrow
(f_a)^{-1}(U)$ such that
\begin{enumerate}
\item $f|_{f^{-1}(U)}=(f_a)|_{f_a^{-1}(U)}\circ h$ and
  $\phi=\phi_a\circ h$, and
\item $(f_a:X_a \rightarrow B, \phi_a)$ is in the same connected component of
  $\Hetc_{S_d}\Betc$ as $(f:X\rightarrow B,\phi)$.
\end{enumerate}
\end{cor}

\begin{proof}
By corollary ~\ref{cor-cor1}, it suffices to consider the case that
$w=2d$.  Suppose the branching monodromy group of $(f:X\rightarrow
B,\phi)$ is $G=S_{A_1}\times \dots \times S_{A_r}$.  We will prove the
result by induction on $r$.  If $r=1$, there is nothing to prove.  So
assume that $r>1$ and assume the result is proved for all smaller
values of $r$.

\

Since $\sum_m (w_m-2\#A_m)$ equals $w-2d=0$, there is some $m$ such
that $w_m\geq 2\#A_m$.  Without loss of generality, suppose $w_1\geq
2\#A_1$.  By proposition~\ref{prop-split}, we may suppose that the
transpositions in $D_1$ are of the form
$(\tau_1,\dots,\tau_{w_1-2},\tau,\tau)$ such that
$\tau_1,\dots,\tau_{w_1-2}$ generate $S_{A_1}$.  But then, choosing a
small disk $D'$ which contains only the branch points $b_{2w_1-1}$ and
$b_{2w_1}$, we see that $(f:X\rightarrow B,\phi)$ and $D'$ satisfy the
hypothesis of lemma~\ref{lem-dbltran}.  Suppose that $j_b\in A_1$ and
$k_b\in A_2$.  By lemma~\ref{lem-dbltran}, we can find a pair
$(f_b:X_b\rightarrow B,\phi_b) \in H^{d,w}\Betc$ and
$h_b:f^{-1}(B-D')\rightarrow (f_b)^{-1}(B-D')$ such that
\begin{enumerate}
\item $f|_{f^{-1}(B-D')} = (f_b)|_{f_b^{-1}(B-D')} \circ h$ and $\phi
  = \phi_b \circ h$, 
\item $\text{br}(f_b)=\text{br}(f)$, 
\item the transposition of $(f_b:X_b\rightarrow B,\phi_b)$ corresponding
  to $\gamma_{w_1-1}$ and $\gamma_{w_1}$ is $(j_b,k_b)$, and
\item $(f_b:X_b \rightarrow B, \phi_b)$ is in the same connected
  component of $\Hetc_{S_d}\Betc$ as $(f:X \rightarrow B, \phi)$.
\end{enumerate}

\

Since the branching monodromy group of $(f_b:X_b \rightarrow B,
\phi_b) \in \Hetc_{S_d}\Betc \cap V$ \emph{outside} of $D'$ already generates
$S_{A_1}\times S_{A_2}$, when we add the transposition $(j_b,k_b)$ we
conclude the
branching monodromy group of $(f_b:X_b \rightarrow B, \phi_b)$ is  
$S_{A_1\cup A_2}\times S_{A_3}\times
\dots \times S_{A_r}$.  By the induction assumption, there is
$(f_a:X_a \rightarrow B,\phi_a)$ and $h_a:f_b^{-1}(U) \rightarrow
f_a^{-1}(U)$ satisfying the conditions of our
 corollary where $(f:X \rightarrow B,\phi)$ is replaced by 
$(f_b:X_b \rightarrow B, \phi_b)$.  Then defining $h=h_a\circ h_b$, we
see that $(f_a:X_a \rightarrow B, \phi_a)$ and $h$ satisfy the
conditions of the corollary for $(f:X \rightarrow B, \phi)$, and the
corollary is proved by induction.
\end{proof}

\section{Induction Argument}~\label{sec-induct}

In this section we will prove that for $w \geq 2d$, $\Hetc_{S_d}\Betc$
is connected.  The basic strategy is as follows:  If $h = g(B) = 0$,
then this is a classical result due to Hurwitz (see the references in
the introduction).  Suppose
given a disk $D\subset B$ and two pairs $(f_1:X_1\rightarrow
B,\phi_1)$ and $(f_2:X_2\rightarrow B, \phi_2)$ such that all branch
points of $f_1$ and $f_2$ are contained in $D$ and such that both
$f_1$ and $f_2$ are trivial over $B-D$, i.e. $f_i^{-1}(B-D)\rightarrow
B-D$ is just $d$ isomorphic copies of $B-D$ for $i=1,2$.  Then the
genus $0$ argument shows that $(f_1:X_1\rightarrow B, \phi_1)$ and
$(f_2:X_2 \rightarrow B, \phi_2)$ are contained in the same connected
component of $\Hetc_{S_d}\Betc$.  So the argument is reduced to
proving that given a general pair $(f:X\rightarrow B, \phi)$ with
branch points in $D$, we can perform braid moves such that
$f^{-1}(B-D)\rightarrow B-D$ is trivial. 

\

\begin{figure}
\begin{picture}(330,210)
\put(10,10){\makebox(330,180){\includegraphics{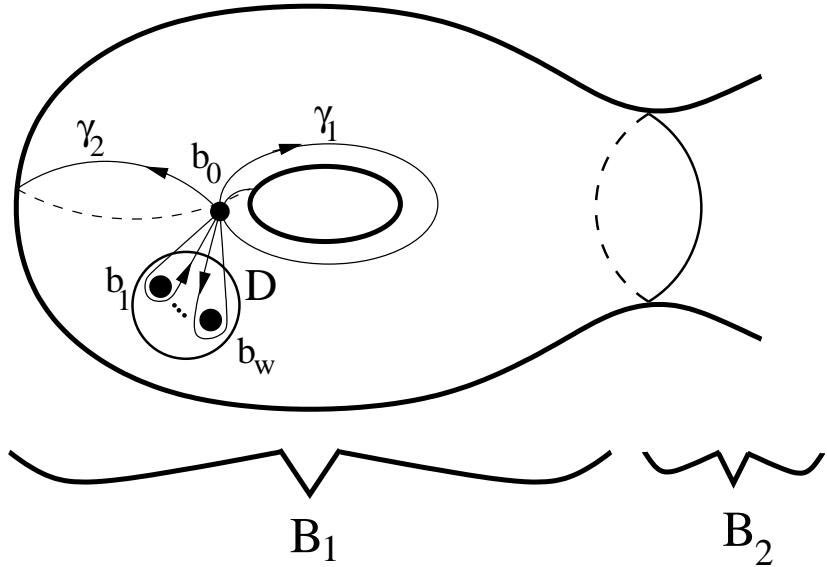}}}
\end{picture}
\caption{The subset $B_1$}
\label{fig-Q1}
\end{figure}

\

Suppose $g \geq 1$ and choose a disk $D\subset B_1\subset B$ situated as in
Figure~\ref{fig-Q1} and which is disjoint from $S$.  Let $V$ be as in
section~\ref{sec-brmon} with respect to this disk $D$.  Every
connected component of $\Hetc_{S_d}\Betc$ clearly intersects $V$.  So
to prove that $\Hetc_{S_d}\Betc$ is connected, it suffices to prove
that for any two pairs $(f_1: X_1\rightarrow B, \phi_1)$ and $(f_2:
X_2\rightarrow B, \phi_2)$ in $\Hetc_{S_d}\Betc \cap V$ 
both pairs in the same connected component.
We prove this by induction on $g$ through
a sequence of intermediate steps (showing each pair is in the same
connected component as a pair with some special properties, and
finally linking up the resulting pairs).

\

\begin{figure}
\begin{picture}(330,210)
\put(10,10){\makebox(330,180){\includegraphics{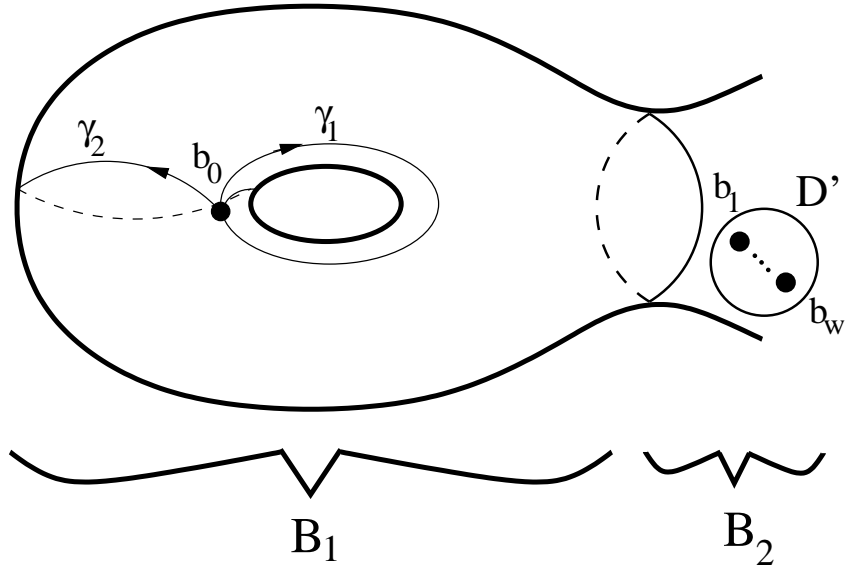}}}
\end{picture}
\caption{The disk $D'$}
\label{fig-Q2}
\end{figure}

\

Let $D'\subset B_2\subset B$ be as in Figure~\ref{fig-Q2}.  Let $V'$ be as in
section~\ref{sec-brmon} with respect to $D'$.  We say that $(f:X
\rightarrow B, \phi)$ in $\Hetc_{S_d}\Betc \cap V'$ is $B_1$-trivial
if $f^{-1}(B_1)\rightarrow B_1$ is a trivial cover.  

\begin{prop}~\label{prop-b1triv}
Suppose $w\geq 2d$.  
Any pair $(f:X\rightarrow B,\phi) \in \Hetc_{S_d}\Betc \cap V$ is in
the same connected component as a pair $(f_a:X_a \rightarrow B,
\phi_a)$ which is $B_1$-trivial.
\end{prop}

\begin{proof}
We prove this in a number of steps.  The idea is to apply braid moves
to reduce $\tilde{\phi}(\gamma_1)$ and $\tilde{\phi}(\gamma_2)$ to the
identity.  Finally we will move all the branch points out of $B_1$
along a specified path to give a $B_1$-trivial pair.

\

\begin{figure}
\begin{picture}(330,210)
\put(10,10){\makebox(330,180){\includegraphics{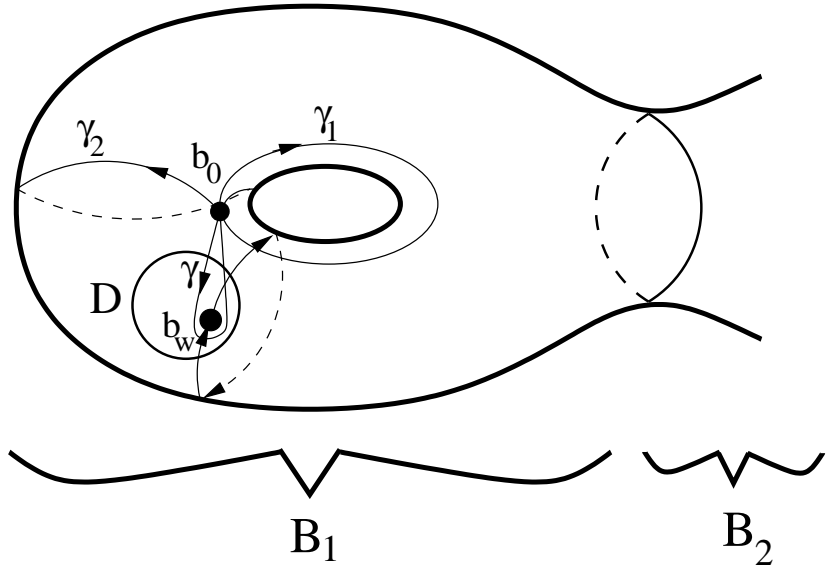}}}
\end{picture}
\caption{The first braid move}
\label{fig-Q3}
\end{figure}
\

Our first braid move is displayed in Figure~\ref{fig-Q3}.  It consists of
choosing the final branch point $b_w$, moving $b_w$
across the loop $\gamma_1$, without crossing
$\gamma_2$, and continuing along the loop ``parallel'' to $\gamma_2$
to return $b_w$ into $D$.  If the resulting cover is $(f_b:X_b
\rightarrow B, \phi_b)$, then we clearly have
$\tilde{\phi_b}(\gamma_1) =
\tilde{\phi}(\gamma_1)\tilde{\phi}(\gamma)$ and
$\tilde{\phi_b}(\gamma_2) = \tilde{\phi}(\gamma_2)$.  So the result is
to multiply the permutation of $\gamma_1$ by the permutation
of $\gamma$ while leaving the permutation of $\gamma_2$ unchanged.

\

\begin{figure}
\begin{picture}(330,210)
\put(10,10){\makebox(330,180){\includegraphics{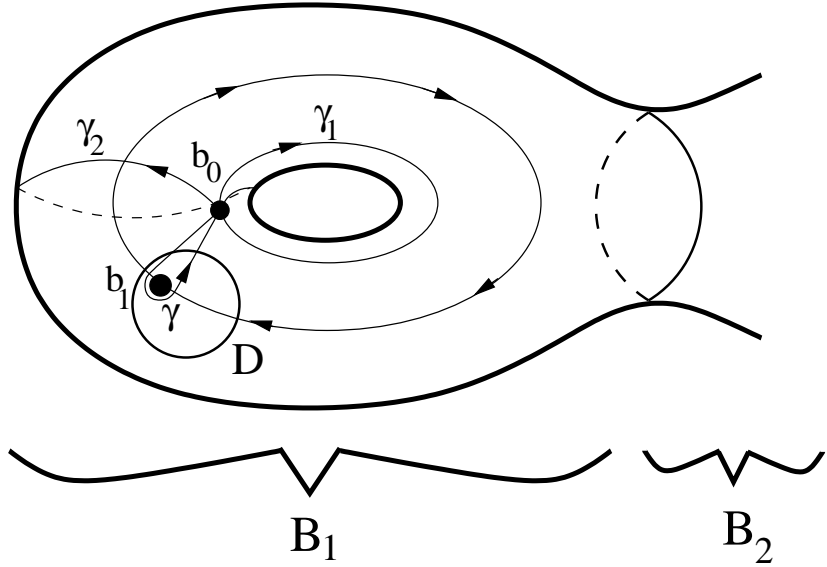}}}
\end{picture}
\caption{The second braid move}
\label{fig-Q4}
\end{figure}
\

Our second braid move is exactly like our first braid move with the
roles of $\gamma_1$ and $\gamma_2$ switched.  
It is illustrated in Figure~\ref{fig-Q4}
We choose the first branch point $b_1$, move $b_1$
across the loop $\gamma_2$, without crossing
$\gamma_1$, and then continue along the loop ``parallel'' to
$\gamma_1$ we return $b_1$ into $D$.  If the resulting cover if
$(f_b:X_b \rightarrow B, \phi_b)$, then we clearly have
$\tilde{\phi_b}(\gamma_2) = \tilde{\phi}(\gamma_2)
\tilde{\phi}(\gamma)$ and $\tilde{\phi_b}(\gamma_1) =
\tilde{\phi}(\gamma_1)$.  So the result is to multiply the permutation
of $\gamma_2$ by the permutation of $\gamma$ while leaving the
permutation of $\gamma_1$ unchanged.  Notice that in both of these
moves, we are not concerned about the effect of the braid move on the
branching monodromy of $D$ (we may always use
corollary~\ref{cor-maincor} to ``repair'' the branching monodromy of $D$).

\

The main claim is that these braid moves along with
corollary~\ref{cor-maincor} suffice to trivialize the permutations of
$\gamma_1$ and $\gamma_2$.  Suppose given $(f:X \rightarrow B, \phi) \in
\Hetc_{S_d}\Betc \cap V$.  Suppose that $\tilde{\phi}(\gamma_1)$ has
cycle type $\lambda = (\lambda_1 \geq \dots \geq \lambda_s)$ and
$\tilde{\phi}(\gamma_2)$ has cycle type $\mu = (\mu_1 \geq \dots \geq
\mu_t)$.  Define $|\lambda| = \sum_m (\lambda_m - 1) = d - s$ and
define $|\mu| = \sum_n (\mu_n - 1) = d - t$.
We claim that there
is a pair $(f_b:X_b \rightarrow B, \phi_b) \in \Hetc_{S_d}\Betc \cap
V$ such that:
\begin{enumerate}
\item $\tilde{\phi_b}(\gamma_1) = \tilde{\phi_b}(\gamma_2) = 1$, and
\item $(f_b:X_b \rightarrow B, \phi_b)$ is contained in the same
  connected component of $\Hetc_{S_d}$ as $(f:X \rightarrow B, \phi)$.
\end{enumerate}
We will prove this by induction on $|\lambda| + |\mu|$.  
If $|\lambda| + |\mu| = 0$, i.e. $\lambda = \mu = 1^d$, we may simply
take $(f_b:X_b \rightarrow B, \phi_b) = (f:X \rightarrow B, \phi)$.
Therefore suppose that $|\lambda| + |\mu| > 0$ and, by way of
induction, suppose the result is proved for all smaller values of
$|\lambda| + |\mu|$.
We make one reduction at the outset:
by corollary~\ref{cor-maincor}, we may replace $(f:X\rightarrow B,
\phi)$ with a pair which is equivalent over $B-D$, but whose branching
monodromy group is all of $S_d$.  

\

Suppose first that $|\lambda| > 0$.  Let
$\sigma \in S_d$ be the $\lambda_1$-cycle occurring in
$\tilde{\phi}(\gamma_1)$ and suppose $\tau \in S_d$ is a transposition
such that $\sigma \tau$ is a $(\lambda_1 - 1)$-cycle.  By
proposition~\ref{prop-split}, we may replace $(f:X\rightarrow B,
\phi)$ by a pair which is equivalent over $B-D$, and whose sequence of
transpositions is of the form $(\tau_1, \dots,
\tau_{w-2},\tau,\tau)$.  
If we apply our first braid move,
the resulting cover $(f_c:X_c \rightarrow B, \phi_c)$ is such that
$|\lambda_c| = |\lambda|-1$ and $|\mu_c| = |\mu|$ so that $|\lambda_c +
|\mu_c| < |\lambda| + |\mu|$.  
By the induction assumption applied to
$(f_c:X_c \rightarrow B, \phi_c)$, we conclude there exists a pair
$(f_b:X_b \rightarrow B, \phi_b)$ with $\tilde{\phi_b}(\gamma_1) =
\tilde{\phi_b}(\gamma_2) = 1$ and which is in the same connected
component of $\Hetc_{S_d}\Betc$ as $(f_c:X_c \rightarrow B, \phi_c)$.
By construction, $(f_c:X_c \rightarrow B, \phi_c)$ is in the same
connected component of $\Hetc_{S_d}\Betc$ as $(f:X\rightarrow B,
\phi)$.  Thus $(f_b:X_b \rightarrow B, \phi_b)$ satisfies conditions
$(1)$ and $(2)$ above.

\

The second possibility is that $|\lambda| = 0$ but $|\mu| > 0$.  Let
$\sigma \in S_d$ be the $\mu_1$-cycle occurring in
$\tilde{\phi}(\gamma_2)$ and suppose $\tau \in S_d$ is a transposition 
such that $\sigma\tau$ is a $(\mu_1-1)$-cycle.  By an obvious
generalization of proposition~\ref{prop-split}, we may replace
$(f:X\rightarrow B, \phi)$ by a pair which is equivalent over $B-D$
and whose sequence of transpositions is of the form
$(\tau,\tau,\tau_1,\dots,\tau_{w-2})$. 
If we apply our
second braid move, the resulting cover $(f_c:X_c \rightarrow B,
\phi_c)$ is such that $|\lambda_c| = |\lambda| = 0$ and $|\mu_c| =
|\mu| - 1$ so that $|\lambda_c| + |\mu_c| < |\lambda| + |\mu|$.  
By the induction assumption applied to
$(f_c:X_c \rightarrow B, \phi_c)$, we conclude there exists a pair
$(f_b:X_b \rightarrow B, \phi_b)$ with $\tilde{\phi_b}(\gamma_1) =
\tilde{\phi_b}(\gamma_2) = 1$ and which is in the same connected
component of $\Hetc_{S_d}\Betc$ as $(f_c:X_c \rightarrow B, \phi_c)$.
By construction, $(f_c:X_c \rightarrow B, \phi_c)$ is in the same
connected component of $\Hetc_{S_d}\Betc$ as $(f:X\rightarrow B,
\phi)$.  Thus $(f_b:X_b \rightarrow B, \phi_b)$ satisfies conditions
$(1)$ and $(2)$ above.  So in both the first and second case, we
conclude that the claim is true for $(f:X \rightarrow B, \phi)$.  So
the claim is proved by induction.

\

\begin{figure}
\begin{picture}(330,210)
\put(10,10){\makebox(330,180){\includegraphics{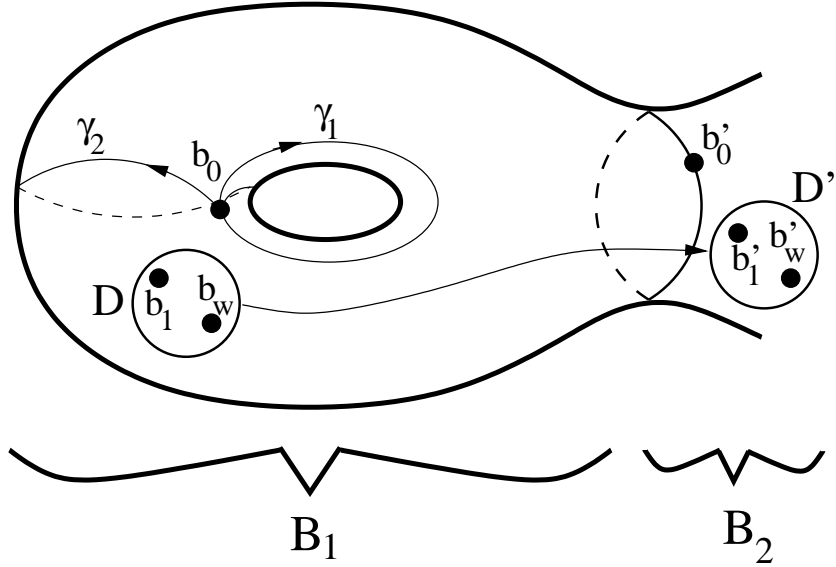}}}
\end{picture}
\caption{Moving $D$ to make the cover $B_1$-trivial}
\label{fig-Q5}
\end{figure}
\

Now we prove the proposition.  By the claim, we may suppose that $(f:X
\rightarrow B, \phi)$ is such that $\tilde{\phi}(\gamma_1) =
\tilde{\phi}(\gamma_2) = 1$.  Finally we move the disk $D$ and all its
branch points out of $B_1$ to $D'$ as shown in Figure~\ref{fig-Q5}.  Let $(f_a:
X_a \rightarrow B, \phi_a)$ be the resulting pair.
Notice that since the path of $D$ never crosses $\gamma_1$ or
$\gamma_2$, we still have $\tilde{\phi_a}(\gamma_1) =
\tilde{\phi_b}(\gamma_2) = 1$.  As the fundamental group
$\pi_1(B_1,b_0)$ is generated by $\gamma_1$ and $\gamma_2$, we
conclude that $(f_a:X_a \rightarrow B, \phi_a)$ is trivial over
$B_1$.  Thus $(f_a:X_a \rightarrow B, \phi_a)$ is $B_1$-trivial, and
the proposition is proved.
\end{proof}

Now we are ready to prove the theorem.

\begin{thm}~\label{thm-main}
If $w \geq 2d$, then $\Hetc_{S_d}\Betc$ is connected.
\end{thm}

\begin{proof}
The proof is by induction on the genus $h$ of $B$.  If $h=0$, the
theorem is due to Hurwitz (see the references in the introduction).  
Thus suppose $h > 0$, and by way of
induction suppose that the theorem is proved for all genera smaller
than $h$.  
Suppose that
$\Sigma\subset \Sigma'\subset B$.  
There is a natural map $\Hetc_{S_d}(B,\Sigma',b_0)
\rightarrow \Hetc_{S_d}\Betc$ whose image is a dense Zariski open
set.  So if $\Hetc_{S_d}(B,\Sigma',b_0)$ is connected, it follows that
$\Hetc_{S_d}\Betc$ is also connected.  Therefore we may enlarge $S$,
if need be, so that it contains a point $b'_0$ in the boundary circle
$B_1\cap B_2$ (and such that this is the only point of $S$ on the
boundary circle).

\

Now by proposition~\ref{prop-b1triv}, we see that every connected
component of $\Hetc_{S_d}\Betc$ contains a $B_1$-trivial pair.  So to
finish the proof, it suffices to prove that for two $B_1$-trivial
pairs, say $(f_1:X_1 \rightarrow B, \phi_1)$ and $(f_2:X_2 \rightarrow
B, \phi_2)$, there are braid moves which change the first pair to the
second.  Let $U\subset B_2$ denote the interior of $B_2$, i.e. the
complement of the boundary circle.  Choose a path $\gamma$ in $B_1$
from $b_0$ to $b'_0$ and in this way identify
$\phi_i:f_i^{-1}(b_0)\rightarrow \{1,\dots,d\}$ with
$\phi'_i:f_i^{-1}(b'_0) \rightarrow \{1,\dots, d\}$.
Now $U$ is homeomorphic to
$B'-\{b'_0\}$ for some Riemann surface $B'$ of genus $h-1$ and for some
point $b'_0 \in B'$.  Let $\Sigma'\subset B'$ denote the union of the
image of $\Sigma\cap U$ and $\{b'_0\}$.  Then the restricted covers
$(f_i: f_i^{-1}(U \rightarrow U, \phi_i)$ for $i=1,2$ are equivalent
to covers $(f'_i:X'_i \rightarrow B',\phi'_i)$ in
$\Hetc_{S_d}(B',\Sigma',b'_0)$.  By the induction assumption, we know that
$\Hetc_{S_d}(B',\Sigma',b'_0)$ is connected.  Therefore there is a path
$\alpha:[0,1] \rightarrow (B'-\Sigma')^0_w$ such that 
\begin{enumerate}
\item $\alpha(0) = \text{br}(f'_1)$,
\item $\alpha(1) = \text{br}(f'_2)$, and
\item if $\tilde{\alpha}:[0,1] \rightarrow \Hetc_{S_d}(B',\Sigma',b'_0)$ is 
the lift with $\tilde{\alpha}(0)=(f'_1:X'_1 \rightarrow B', \phi'_1)$,
then $\tilde{\alpha}(1) = (f'_2:X'_1 \rightarrow B', \phi'_2)$. 
\end{enumerate}

\

Using our homeomorphism, we may identify $\alpha$ with a path
$\beta:[0,1] \rightarrow (U-\Sigma\cap U)^0_w$ such that $\beta(0) =
\text{br}(f_1)$ and $\beta(1)=\text{br}(f_2)$.  It follows that if
$\tilde{\beta}:[0,1] \rightarrow \Hetc_{S_d}\Betc$ is the lift with
$\tilde{\beta}(0) = (f_1:X_1 \rightarrow B, \phi_1)$, then
$\tilde{\beta}(1) = (f_2:X_2 \rightarrow B, \phi_2)$.  This proves
that $(f_1:X_1 \rightarrow B, \phi_1)$ and $(f_2:X_2 \rightarrow B,
\phi_2)$ lie in the same connected component of $\Hetc_{S_d}\Betc$.
It follows that $\Hetc_{S_d}\Betc$ is connected, and the theorem is
proved by induction.
\end{proof}

Now we can prove theorem~\ref{thm-thm1}.  There is a forgetful map
$\Hetc_{S_d}\Betc \rightarrow {\mc H}^{d,w}_{S_d}(B)$.  This morphism
is \'etale with dense image.  Since $\Hetc_{S_d}\Betc$ is connected,
it follows that ${\mc H}^{d,w}_{S_d}(B)$ is also connected, which
proves theorem~\ref{thm-thm1}.

\bibliography{my}
\bibliographystyle{abbrv}

\end{document}